\documentclass[a4paper,twoside]{article}
\usepackage{amsfonts,amssymb,mathrsfs,amsmath}
\usepackage[dvips]{color}
\usepackage{graphicx}
\title{\sc Rigidity and the Chess Board Theorem for Cube Packings}
\date{}
\author{Andrzej P. Kisielewicz
and Krzysztof Przes{\l}awski\\
\\
{\small Wydzia{\l} Matematyki, Informatyki i Ekonometrii, Uniwersytet Zielonog\'orski}\\
{\small ul. Z. Szafrana 4a, 65-516 Zielona G\'ora, Poland}\\
{\small A.Kisielewicz@wmie.uz.zgora.pl}\\
{\small K.Przeslawski@wmie.uz.zgora.pl}}
\newtheorem{pr}{\sc Proposition}

\newtheorem{tw}[pr]{\sc Theorem}
\newtheorem{wn}[pr]{\sc Corollary}
\newtheorem{df}{\sc Definition}

\newtheorem{uw}{\sc Remark}
\newtheorem{uwi}[uw]{\sc Remarks}
\newenvironment{uwa}{\begin{uw} \rm}{\end{uw}}

\newtheorem{nap}{\sc Example}

\newtheorem{nps}[nap]{\rm EXAMPLES}

\def\ka #1{\mathscr{#1}}
\def\kal #1 #2{\mathscr{#1}^{#2}}
\def\proof{\noindent \textit{Proof.\,\,\,}}
\def\dowod #1{\noindent\textit{Proof #1.\,\,\,}}

\def\zet{\mathbb{Z}}

\def\er{\mathbb{R}}

\def\ind #1{\operatorname{ind}(#1)}
\def\sgn #1{\operatorname{sgn}(#1)}

\begin{document}
\maketitle
\begin{abstract}
Each packing of $\er^d$ by translates of the unit cube $[0,1)^d$ admits a decomposition   
into at most two parts such that if  a translate of the unit cube is covered by one of them, then it also belongs to such a part. 

\end{abstract}

\medskip
%A non-empty  subset $T$ of $\er^d$ is said to be a \textit{Keller set} if for every two vectors $s$, $t\in T$ there is %$i\in [d]$ such that $s_i-t_i\in \zet\setminus\{0\}$. 

%Let $I=[0,1)^d$. Observe that the fact that $T$ is a Keller set implies that the set $I+T = \{I+t: t\in T\}$ consists of disjoint sets, that is, it is a packing of $\er^d$ by half-open unit cubes. It is also a cube tiling of $\bigcup (I+T)$. 
%We shall  refer to $\bigcup I+T$ as a \textit{Keller polybox}. 

Let  $I=[0,1)^d$, and  $S\subset \er^d$ be a  non-empty set. We say that the set $I+S = \{I+s\colon s\in S\}$ is a \textit{packing}  of $\er^d$  by (half-open) unit cubes if the members of $I+S$ are pairwise disjoint. Clearly, $I+S$ is a packing if and only  if for every two vectors  $t$, $t'\in S$ there is $i\in [d]$ such that $|t_i-t'_i|\ge 1$.
The cube-packing $I+S$ of $\er^d$ is a \textit{tiling} of $F\subseteq\er^d$ if $F=\bigcup (I+S)$. The set $F$ is said to be \textit{rough} if for each $u\in \er^d$ the inclusion  $I+u\subseteq F$ implies that $u\in S$. Moreover, $F$ is \textit{rigid} if $I+S$ is a unique tiling of $F$. 

Lagarias and Shor conjectured in \cite{LS} that if  $I+T$ is a \textit{2-extremal cube-tiling} of $\er^d$, then $T$ decomposes into two explicitly defined parts $T^0$ and $T^1$ such that each of them determines $T$. 
Clearly, the conjecture means simply that the sets $F^0:=\bigcup (I+T^0)$ and $F^1:=\bigcup(I+T^1)$ are rigid in the sens defined previously. In \cite{KP}, we confirmed the conjecture proving a more  general result on the rigidity of \textit{polyboxes}. (We believe that the notion of a polybox appeared therein for the first time.)  One of the referees of our work acting for the \textit{Discrete and Computational Geometry} asked us whether there is a more straightforward approach to the rigidity of tilings. We hope that the present work answers this question in 
the positive. We offer two proofs of our main result (Theorem \ref{rt}). Interestingly enough, they resemble two of the fourteen proofs (collected by S. Wagon \cite{W}) of the result that if a rectangle can be tiled by rectangles each of which has at least one integer side, then the tiled rectangle has an integer side. Our first proof resembles that of A. Douady, while the other,  that of R. Rochberg, and S. K. Stein. The conjecture of Lagarias and Shor is an immediate consequence of Theorem \ref{holenderka} which is formulated for packings rather than for tilings. Since their work relates to Keller's conjecture on cube tilings, we give applications of our results to this sort of problems.

\begin{tw}[rigidity theorem]
\label{rt}
Let $I+S$ be a packing of $\er^d$. Suppose that for every two vectors $t$, $t'\in S$ if $t-t'\in \{-1,0,1\}^d$, then 
the number $|\{i\colon |t_i-t'_i|=1\}|$ is even. Then $F=\bigcup(I+S)$ is rough. In particular, $F$ is rigid.
\end{tw}

Clearly, the rigidity theorem can be rephrased as follows:
\begin{tw}
\label{pokrycie}
Let $I+S$ be a packing of $\er^d$  and  let $u\not\in S$. If  $I+u \subseteq\bigcup(I+S)$, then there are $t$, $t'\in S$ such that 
$t-t'\in \{-1,0,1\}^d$ and the number $|\{i\colon |t_i-t'_i|=1\}|$ is odd.
\end{tw}

\dowod 1 We let $u=0$ without loss of generality. In addition, we may assume that each of the cubes $I+s$, $s\in S$,  intersects $I$. For every $i\in [d]$, let 
$$
A_i=\{a>0\colon \text{there is $s\in S$ such that $a=s_i$}\}
$$
and
$$
B_i=\{b<0\colon \text{there is $s\in S$ such that $b=s_i$}\}.
$$  
Let us remark that if $b\in B_i$, then $b+1\in A_i$. In order to show this fact, assume for simplicity that $i=d$. Fix any $v\in S$ such that $v_d=b$. Pick $y\in I\cap (I+s)$. Let $z=(y_1,\ldots , y_{d-1})$. Clearly, the line segment $J=\{z\}\times [0,1)$ intersects with at most two of the sets belonging to $I+S$. Since $I+v$
is one of them and $J\not\subseteq I+v$, there is yet another one, say $I+w$. These two cubes cover $J$. Since
$b<0$, we get $J\cap(I+v)=\{z\}\times [0,b+1)$. Consequently, $J\cap(I+w)=\{z\}\times [b+1, 1)$, which implies $w_i=b+1\in A_i$.

Let $U=\{1\}\times A_1\cup\cdots\cup \{d\}\times A_d$. 
%For each $(i,a)\in U$ and each $x\in\er^U$, let us write $x_{ia}$ instead of $x(i,a)$.  
Then, by the above remark, for every $s\in S$, we can define the mapping $f_s\colon \er^U\to\er^d$ by the formula
$$
(f_s(x))_i=\left\{
\begin{array}{ll}
x(i,s_i+1)-1 & \text{if $s_i<0$,}\\
0 & \text{if $s_i=0$,}\\
x(i,s_i) & \text{if $s_i>0$}. 
\end{array}
\right.
$$
Let us define $p\in \er^U$ so that $p(i,a)=a$, $(i,a)\in  U$. Observe that $f_s(p)=s$ for all $s\in S$. Let us define $\varepsilon>0$ so that for each $i$ and each $a\in A_i$, the interval $(a - \varepsilon ,a +\varepsilon)$ does not contain any element of the set $\{0,1\}\cup (A_i\setminus \{a\})$.
Let 
$$
V=\{q\in \er^U\colon \text{$|q(i,a)-p(i,a)|<\varepsilon $, for every $(i,a)\in U$}\}.
$$
Let $W$ be the set of all these $q\in V$ for which $I+f_s(q)$, $s\in S$, are disjoint cubes intersecting  $I$, whose union contains $I$. As $p\in W$, it follows that $W$ is non-empty. We prove now that $V=W$.

Let $q^0\in W$. Fix $(i,a)\in U$ and pick any element $q^1\in V$ such that $q^0(j,c)=q^1(j,c)$, whenever $(j,c)\neq(i,a)$. Let $R^\tau=\{f_s(q^\tau)\colon s\in S\}$, where $\tau\in\{0,1\}$. Define 
$$
R^\tau_a=\{f_s(q^\tau)\colon (f_s(q^\tau))_i=q^\tau(i,a)\}.
$$
By the definition of $\varepsilon$, $q^\tau(i,a)\in (0,1)$. For the sake of brevity, assume that $i=d$. Therefore, since $I+R^0$ is a packing of $\er^d$ and at the same time a covering of $I$, it is easily observed that there is a set $X\subseteq [0,1)^{d-1}$ such that 
$$
I\cap\bigcup (I+R^0_a)=X\times [0,1)\,\, \text{and}\,\,  I\cap\bigcup (I+(R^0\setminus R^0_a))=([0,1)^{d-1}\setminus X)\times [0,1).
$$  
Again, by the definition of $\varepsilon$,
$$
R^1_a=R^0_a+ (q^1(i,a)-q^0(i,a))e_i,
$$
where $e_i$ is the $i$-th vector of the standard basis of $\er^d$. This equality together with the already mentioned relation $q^\tau(i,a)\in (0,1)$ imply that $I\cap\bigcup (I+R^1_a)=X\times [0,1)$. Moreover, by the definition of $q^\tau$, the set $R\setminus R^\tau_a$ is independent of $\tau$. Consequently, $I+R^1$ is a covering of $I$, and a packing of $\er^d$. Equivalently, $q^1\in W$. Since $V$ is a cube in $\er^U$, what we have shown implies $V=W$. 
    
Let $x\in V$, and $s\in S$. It follows from the definition of $\varepsilon$ and $f_s$ that $s_i$ and $(f_s(x))_i$ have the same sign. Let $\operatorname{vol}$ denote the standard volume measure in $\er^d$. Therefore, we have
$$
\operatorname{vol}(I\cap (I+f_s(x)))=\prod_{i\colon s_i>0}(1-x(i,s_i))\cdot\prod_{i\colon s_i<0}x(i,s_i+1)=:P_s(x). 
$$ 
Summing up with respect to $s$ gives us
$$
1=\sum_{s}P_s(x)=:P(x).
$$
Thus, $P$ is a polynomial on $\er^d$ which is constant on $V$. Since the latter set is open,  $P$ is constant.    
We have assumed that $0\not\in S$. This implies that  each of the polynomials $P_s$, $s\in S$, is of a positive degree. It is clear that in the expansion  of $P_s$ into monomials there is only one leading term, that is the term of the greatest degree, 
$$
Q_s(x)=(-1)^{|\{i\colon s_i>0\}|}\prod_{i\colon s_i>0}x(i,s_i)\cdot\prod_{i\colon s_i<0}x(i,s_i+1).
$$   
Let us pick $t\in S$ so that $P_t$ is of maximal degree. Since $P$ is constant, we deduce that there is  $t'\in S\setminus\{t\}$ such that $Q_t+Q_{t'}=0$. It is easily seen that this equation is equivalent to saying that $t-t'\in \{-1,0,1\}^d$ and the number $|\{i\colon |t_i-t'_i|=1\}|$ is odd. 

\medskip
\dowod 2
As before, we assume that $u=0$ and $I+s$ intersects $I$, whenever $s\in S$. Choose $t\in S$ so that it has the maximum number of non-zero coordinates. As we can change the order of coordinates if necessary, we may assume that there is $k\in [d]$ such that $t_i\neq 0$, if $i\le k$, otherwise $t_i=0$. Let us identify $\er^k$ with the subspace $\er^k\times \{0\}^{d-k}$. For each $x\in \er^d$, let $x^\dagger=(x_1,\ldots, x_k)$ be the projection of $x$ onto $\er^k$.    
Let us define $v\in \er^k$ as follows 
$$
v_i=\left\{
\begin{array}{ll}
t_i+1, & \text{if $t_i<0$;}\\
t_i, & \text{otherwise}.
\end{array}
\right.
$$
Clearly, $v$ is an interior point of the $k$-dimensional unit cube $I^\dagger = [0,1)^k$. For sufficiently small $\varepsilon>0$, the cube  $\varepsilon I^\dagger + v$ is contained in $I^\dagger$. 
In particular, it is covered by the cubes $I^\dagger + s^\dagger$, for $s\in S$. Let us define 
$$
T=\{s^\dagger\colon \text{$s\in S$,\, $(I+s)\cap (\varepsilon I^\dagger + v)\neq \emptyset$,\, for every $\varepsilon>0$}\}.   
$$
(We identify here $I^\dagger$ with $I^\dagger\times\{0\}^{d-k}$ and we interpret $v$ as an element of $\er^d$ according to the convention we have made.)
It is clear that $I^\dagger+T$ is a packing of $\er^k$ by unit cubes. Moreover, there is $\gamma>0$ such that $B=\gamma I^\dagger+v$ is covered by $I^\dagger+T$, and, at the same time,  included in $I^\dagger$. Let us split  each factor $B_i$ of the cube $B$ into segments $B_i^{0}=[-\gamma,0)+v_i$ and $B_i^1=[0,\gamma)+v_i$. Then $B$ decomposes into $2^n$ cubes $B^\sigma=B^{\sigma_1}_1\times\cdots\times B^{\sigma_d}_d$, where $\sigma\in \{0,1\}^d$ . Let us subordinate to each $B^\sigma$ its \textit{sign} $\sgn{B^\sigma}=(-1)^{\sum_i\sigma_i}$. For every set $C\subseteq B$ that can be represented as a union of a non-empty family $\ka F\subseteq\ka B:=\{B^\sigma\colon \sigma\in\{0,1\}^d\}$, we define the \textit{index} of $C$
$\ind{C}=\sum_{B^\sigma\subseteq C}\sgn{B^\sigma}$. Observe 
that if $C=(I^\dagger+z)\cap B$, where $z\in T$, then $C$ is a box and  for each $i\le k$, the factor $C_i$  is equal to one of the three sets $B_i$, $B_i^0$, $B_i^1$. In  particular, the index of $C$ is well-defined, and is equal to zero if and only if there is $i$ such that $C_i=B_i$, which in turn is equivalent to saying that $z_i=0$. Furthermore, if $C_i\neq B_i$  for each $i$, then there is $\kappa\in \{0,1\}^k$ such that $C=B^\kappa$ and 
$$
\ind {C}=\sgn{B^\kappa}=(-1)^{|\{i\colon z_i=t_i\}|}= (-1)^{|\{i\colon z_i>0\}|},
$$     
as by the definition of $T$, the sets $\{i\colon z_i=t_i\}$ and $\{i\colon z_i>0\}$ have to be equal. Let $T'$ be the subset of $T$ consisting of all $z\in T$ such that $z_i\neq 0$, whenever $i\in [k]$. We have
$$
0=\ind{B}=\sum_{z\in T} \ind{(I^\dagger+z)\cap B}= \sum_{z\in T'} \ind{(I^\dagger+z)\cap B}=\sum_{z\in T'} (-1)^{|\{i\colon z_i>0\}|}.
$$
Since $t=t^\dagger$ belongs to $T'$, there has to exists an additional element $w$ belonging to this set such that the sets $\{i\colon t_i>0\}$ and $\{i\colon w_i>0\}$ are of different parity. Let $t'$ be the element of $S$ which projects on $w$. Since $w_i\neq 0$ for each $i\le k$, then by the definition of $k$, we deduce that $t'=w$. It is easily seen that just constructed $t$ and $t'$ satisfy our conclusion.   

The main idea of the proof can be easily grasped by analyzing the picture below. 

{\center
\includegraphics[width=4cm]{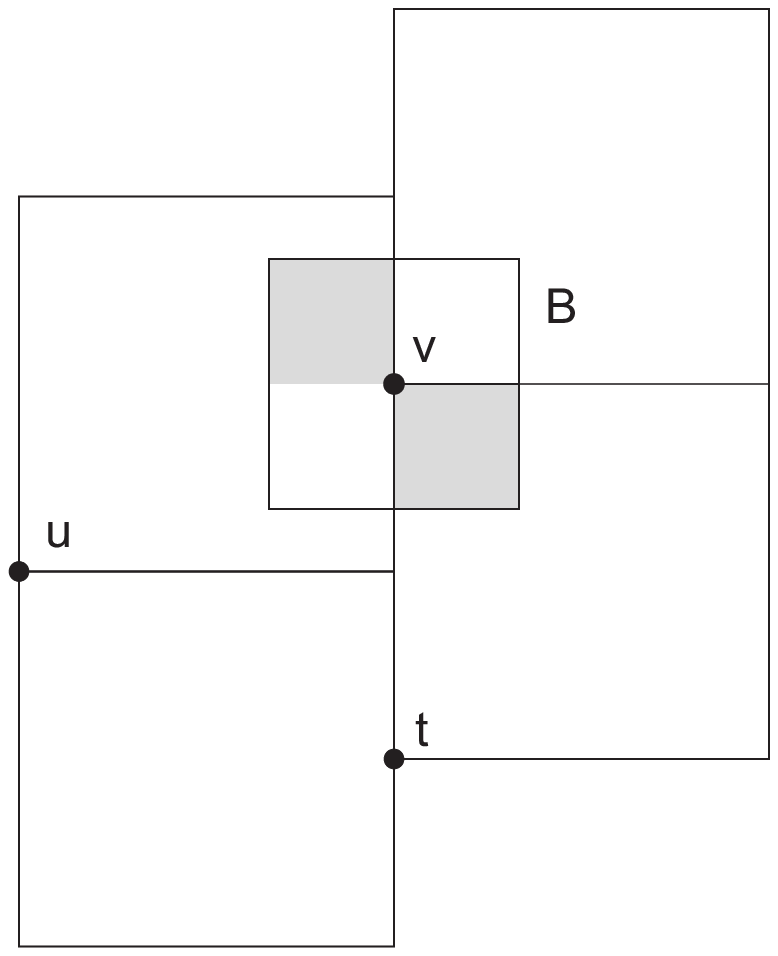}\\
}

\medskip
\noindent{\footnotesize Fig. 1. This picture corresponds to the case $d=k=2$. The shaded areas indicate all these cubes $B^{\sigma}$ in the decomposition of $B$ for which $\sgn {B^\sigma}=-1$. The set $T$ is equal here to $\{t,u,v\}$. Clearly, $\ind{(I+t)\cap B}=-1$, $\ind{(I+u)\cap B}=0$. The element $t'$, which existence is guaranteed by our reasoning, coincides with $v$.   
  }
\hfill $\square$  
\begin{uwa}
\label{silnawersja}
An inspection of any of the two proofs reveals that the conclusion of Theorem \ref{pokrycie} can be strengthen: 
Let $S'$ be the set of all these $s\in S$ for which the intersection $(I+s)\cap (I+u)$ is non-empty , and  $\ka J$   be the family of  all sets 
$\langle s\rangle =\{i\in[d]\colon 0<|s_i-u_i|\}$, $s\in S'$. Then for each $t\in S'$ such that $\langle t\rangle$ is a maximal element of $\ka J$ with respect to the partial order defined by the inclusion, there is $t'\in S'$ such that $\langle t'\rangle=\langle t\rangle$, $t-t'\in \{-1,0,1\}^d$ and $\{i\colon |t_i-t'_i|=1\}$ is of odd cardinality.
\end{uwa}
\begin{tw}[chess board theorem]
\label{holenderka}
If $I+S$ is a cube packing of $\er^d$, then there is a decomposition $S^0$, $S^1$ of $S$  such that the sets $F^i=\bigcup(I+S^i)$, are rough. The sets $S^i$ can be defined explicitly.
\end{tw}      

\proof Let us define two relations $\sim$ and  $\approx$ in $S$ as follows: 
$$
\text{$s\sim t$ if and only if  $s-t\in \zet^d$,} 
$$
$$
\text{$s\approx t$ if and only if $s\sim t$ and the number $|\{i\colon t_i-s_i\equiv 1\, \operatorname{mod} 2 \}|$ is even.}  
$$ 
Both of these relations are equivalences. Either each equivalence class $A$ of the relation $\sim$ is an equivalence class  of $\approx$ or it splits into exactly  two such classes $A'$ and $A''$. Let us pick $S^0$ so that if $A$ does not split, then $A$ is contained in $S^0$ or is disjoint with this set; otherwise, $S^0$ includes exactly one of the classes $A'$, $A''$. It is easily observed that the sets $S^0$ and $S^1=S\setminus S^0$ satisfy the assumptions of the rigidity theorem, therefore they define the desired decomposition.  \hfill$\square$

{\center
\includegraphics[width=5cm]{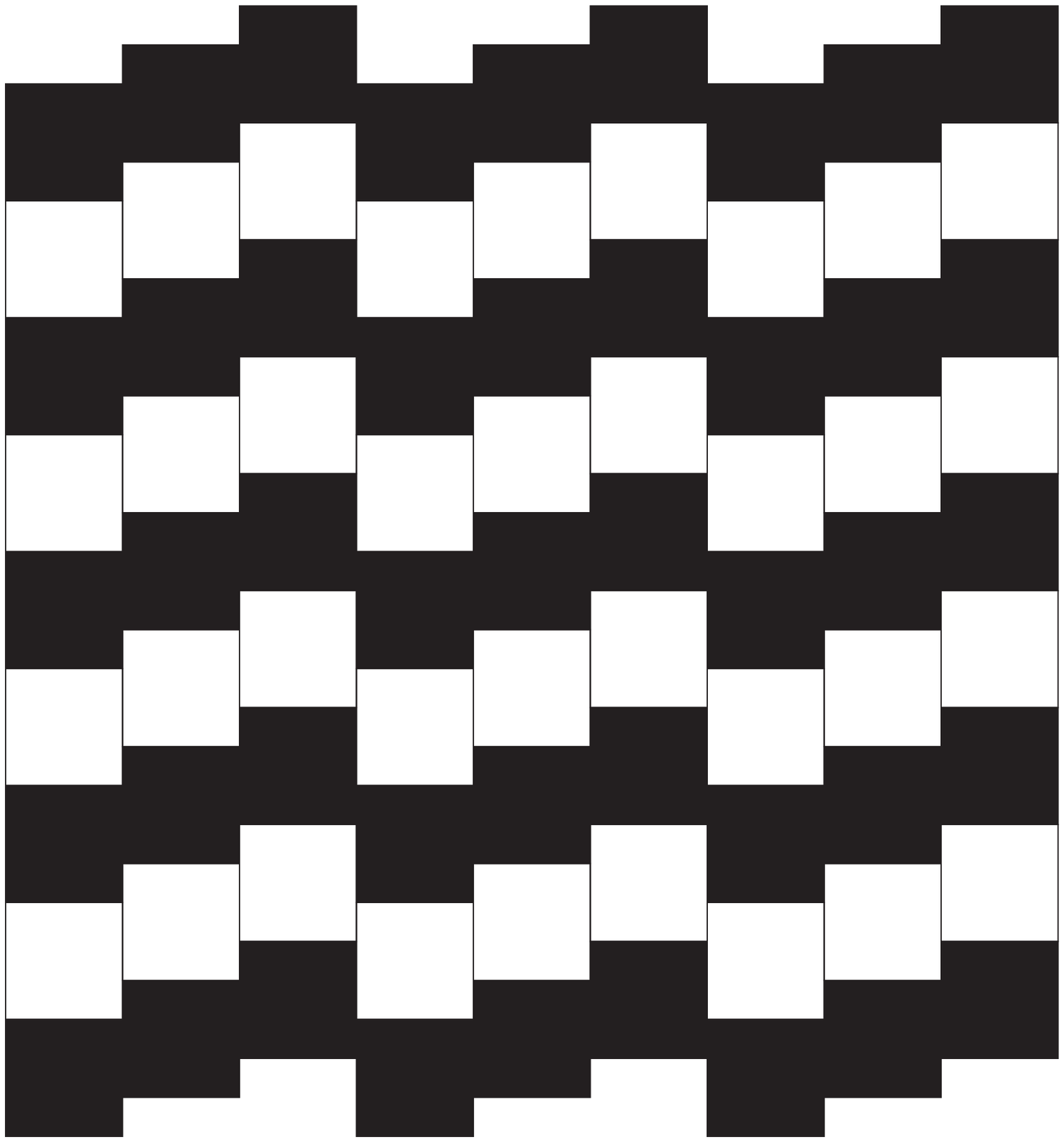}\\
}
\medskip
\noindent{\footnotesize Fig. 2. A cube tiling  and a related decomposition of $\er^2$ into two rough parts: $F^0$ (white) and $F^1$  (black).  
  }

\medskip
Let us remark that one of the sets $S^i$ is allowed to be empty.  

The stronger version of Theorem \ref{pokrycie} described  in Remark \ref{silnawersja} can be applied to cube tilings of $\er^d$. 

\begin{tw}
\label{orthant}
Suppose $I+S$ is a cube tiling of $\er^d$. For every $t\in S$ and  every $\varepsilon\in\{-1, 1\}^d$ there is  a set $J\subseteq [d]$ of odd cardinality such that the vector $t'=t+\sum_{i\in J} \varepsilon_ie_i$, where $e_i$ are elements of the standard basis $e_1=(1,0,\ldots, 0),\ldots, e_d=(0,\ldots, 1)$, belongs to $S$.     
\end{tw}
\proof 
Clearly, we may assume $t=0$. Let  $u=(1/2)\varepsilon$.  Since $I+S$ is a tiling, $I+u$ is contained in $I+S$. Define $S'$ and $\ka J$ as in Remark 1. It is easily seen that $0\in S'$ and $\langle 0\rangle=[d]$. Consequently, the latter set is maximal in $\ka J$. Thus, by Remark 1, there is $t'\in S'$  such that $\langle t'\rangle=[d]$,  and the set 
$J=\{i\in [d]\colon |t'_i-t_i|=|t'_i|=1\}$ is odd. Observe that $\varepsilon_i=-1$, for $i\in J$ if and only if  $t_i'=-1$, as in other case $I+t'$ and $I+u$ would be disjoint. Thus, $t'=\sum_{i\in J} \varepsilon_ie_i$. \hfill$\square$              

\begin{wn}
If $I+S$ is a cube tiling of $\er^d$, then for every $t\in S$ the set $(t+\zet^d)\cap S$ is infinite.
\end{wn} 
 
\pagebreak
{\center
\includegraphics[width=5cm]{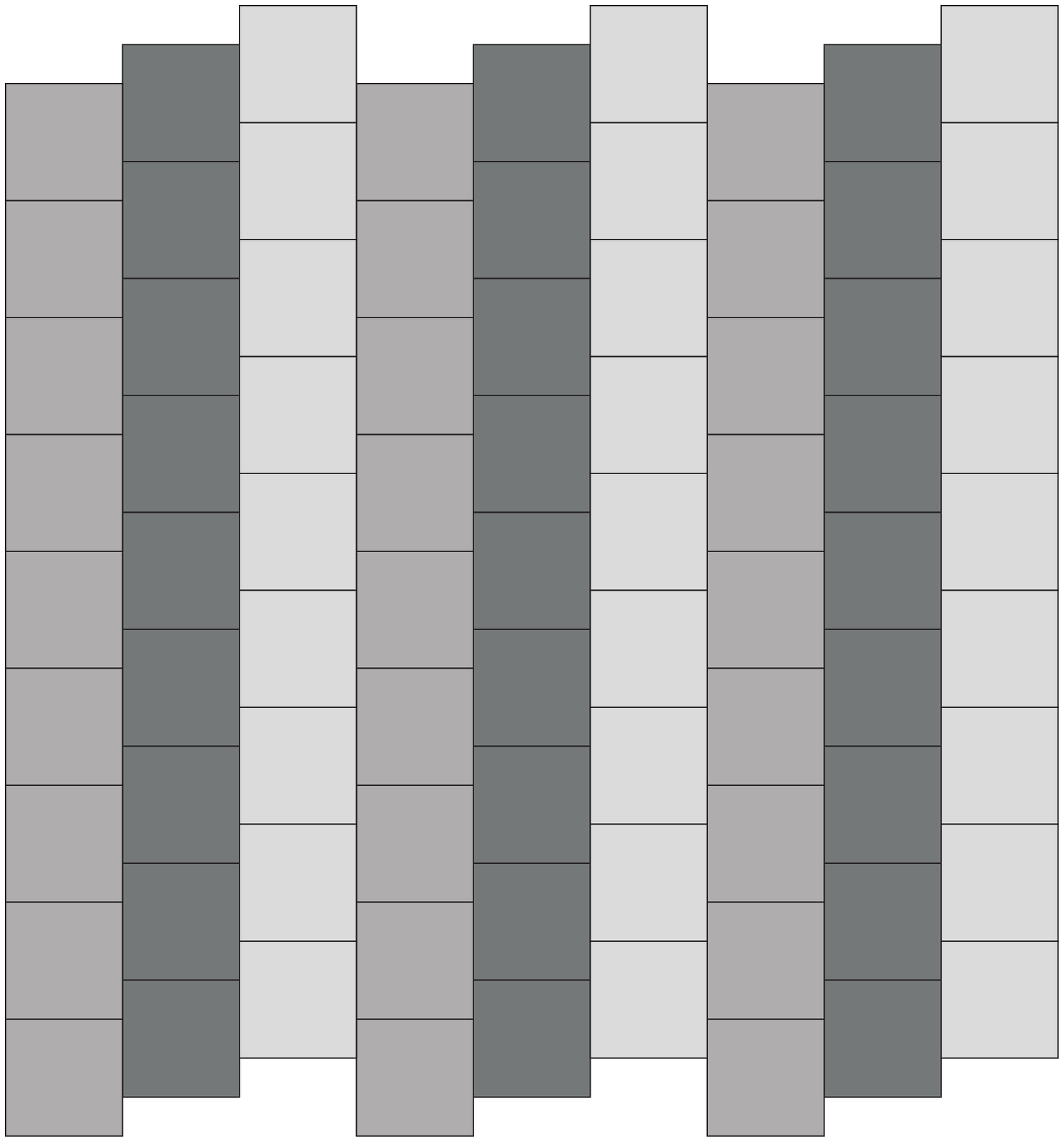}\\
}
\medskip

\noindent{\footnotesize Fig. 3. Each shade of gray represents one of the family of boxes $I+(t+\zet^d)\cap S$, $t\in S$. 
  }

\begin{pr} 
Suppose $I+S$ is a cube tiling of $\er^d$. If $G=S\cap \zet^d$ is a subgroup of $\zet^d$ and there are $k=(k_1,\ldots, k_d)\in \zet^d$  and  a  set $L\subseteq [d]$ containing at least $d-2$ elements  such that 
\begin{itemize}
\item[(1)]for every $i\in [d]$, the multiple $k_ie_i$ of $e_i$ belongs to $G$,
\item[(2)] for every   $i\in [d]$ and $l\in L$, the coordinates  $k_i$ and $k_l$ are relatively prime, whenever $i\neq l$, 
\end{itemize}
then there is $m\in [d]$ such that $e_m\in  G$.        
\end{pr}

\proof
By the preceding theorem and the fact that $0\in S$, there is  a set $J\subseteq [d]$ of odd cardinality such that $s=\sum_{i\in J} e_i$ belongs to $S$. Obviously, $s$ is also an element of $G$. If $J$ is a singleton, then $s$ is a vector of the standard basis; therefore, it remains to consider the case $|J|\ge 3$. Then there is $m\in J\cap L$. Let $n=\prod_{i\in J\setminus\{m\}}k_i$. It follows from assumption (2) that $k_m$ and $n$ are relatively prime. Thus, there exist nonzero integers $x$ and $y$ such that $xn+yk_m=1$. We have $ns= ne_m+\sum_{i\in J\setminus\{m\}}ne_i$. Since the elements $ne_i$, $i\in J\setminus\{m\}$, are multiples of $k_ie_i$, they belong to $G$. Consequently, $ne_m$ belongs to $G$. Now, we have $e_m=x(ne_m)+y(k_me_m)$ belongs to $G$.\hfill$\square$          

\medskip
Theorem \ref{holenderka} is a generalization of Theorem 50 in \cite{KP}. It should  be mentioned however that it can be proved within the framework of a theory developed there. Theorem \ref{pokrycie} relates to Lemma 31 in \cite{KP}. These  results rest upon an idea  which has been already exploited in \cite{BFF} (see also \cite{Z}).

\medskip  
\noindent\textbf{Acknowledgements.} We wish to thank Ron Holzman for pointing out papers of A. Berger, A. Felzenbaum and A. Fraenkel on covering systems to us.

\end{document}